# A new approach to inverse spectral theory, I. Fundamental formalism

By BARRY SIMON

**Abstract**

We present a new approach (distinct from Gel′fand-Levitan) to the theorem of Borg-Marchenko that the $m$-function (equivalently, spectral measure) for a finite interval or half-line Schrödinger operator determines the potential. Our approach is an analog of the continued fraction approach for the moment problem. We prove there is a representation for the $m$-function $m(-\kappa^2) = -\kappa - \int_0^b A(\alpha)e^{-2\alpha\kappa}\,d\alpha + O(e^{-(2b-\varepsilon)\kappa})$. $A$ on $[0,a]$ is a function of $q$ on $[0,a]$ and vice-versa. A key role is played by a differential equation that $A$ obeys after allowing $x$-dependence:
$$\frac{\partial A}{\partial x} = \frac{\partial A}{\partial \alpha} + \int_0^\alpha A(\beta, x) A(\alpha - \beta, x)\, d\beta.$$
Among our new results are necessary and sufficient conditions on the $m$-functions for potentials $q_1$ and $q_2$ for $q_1$ to equal $q_2$ on $[0,a]$.

## 1. Introduction

Inverse spectral methods have been actively studied in the past years both via their relevance in a variety of applications and their connection to the KdV equation. A major role is played by the Gel′fand-Levitan equations. Our goal in this paper is to present a new approach to their basic results that we expect will lead to resolution of some of the remaining open questions in one-dimensional inverse spectral theory. We will introduce a new basic object (see (1.24) below), the remarkable equation, (1.28), it obeys and illustrate with several new results.

To present these new results, we will first describe the problems we discuss. We will consider differential operators on either $L^2(0,b)$ with $b < \infty$ or $L^2(0,\infty)$ of the form

(1.1) $$-\frac{d^2}{dx^2} + q(x).$$

If $b$ is finite, we suppose

$$\beta_1 \equiv \int_0^b |q(x)|\,dx < \infty \tag{1.2}$$

and place a boundary condition

$$u'(b) + hu(b) = 0, \tag{1.3}$$

where $h \in \mathbb{R} \cup \{\infty\}$ with $h = \infty$ shorthand for the Dirichlet condition $u(b) = 0$. If $b = \infty$, we suppose

$$\int_y^{y+1} |q(x)|\,dx < \infty \quad \text{for all } y \tag{1.4}$$

and

$$\beta_2 \equiv \sup_{y>0} \int_y^{y+1} \max(q(x),0)\,dx < \infty. \tag{1.5}$$

Under condition (1.5), it is known that (1.1) is the limit point at infinity [15].

In either case, for each $z \in \mathbb{C}\setminus[\beta, \infty)$ with $-\beta$ sufficiently large, there is a unique solution (up to an overall constant), $u(x, z)$, of $-u'' + qu = zu$ which obeys (1.3) at $b$ if $b < \infty$ or which is $L^2$ at $\infty$ if $b = \infty$. The principal $m$-function $m(z)$ is defined by

$$m(z) = \frac{u'(0, z)}{u(0, z)}. \tag{1.6}$$

We will sometimes need to indicate the $q$-dependence explicitly and write $m(z; q)$. If $b < \infty$, "$q$" is intended to include all of $q$ on $(0, b)$, $b$, and the value of $h$.

If we replace $b$ by $b_1 = b - x_0$ with $x_0 \in (0, b)$ and let $q(s) = q(x_0 + s)$ for $s \in (0, b_1)$, we get a new $m$-function we will denote by $m(z, x_0)$. It is given by

$$m(z, x) = \frac{u'(x, z)}{u(x, z)}. \tag{1.7}$$

$m(z, x)$ obeys the Riccati equation

$$\frac{dm}{dx} = q(x) - z - m^2(z, x). \tag{1.8}$$

Obviously, $m(z, x)$ only depends on $q$ on $(x, b)$ (and on $h$ if $b < \infty$). A basic result of the inverse theory says that the converse is true:

THEOREM 1.1 (Borg [3], Marchenko [12]). *$m$ determines $q$. Explicitly, if $q_1, q_2$ are two potentials and $m_1(z) = m_2(z)$, then $q_1 \equiv q_2$ (including $h_1 = h_2$).*

We will improve this as follows:

THEOREM 1.2. *If $(q_1, b_1, h_1)$, $(q_2, b_1, h_2)$ are two potentials and $a < \min(b_1, b_2)$, and if*

$$(1.9) \qquad q_1(x) = q_2(x) \quad \text{on } (0, a),$$

*then as $\kappa \to \infty$,*

$$(1.10) \qquad m_1(-\kappa^2) - m_2(-\kappa^2) = \tilde{O}(e^{-2\kappa a}).$$

*Conversely, if (1.10) holds, then (1.9) holds.*

In (1.10), we use the symbol $\tilde{O}$ defined by $f = \tilde{O}(g)$ as $x \to x_0$ (where $\lim_{x \to x_0} g(x) = 0$) if and only if $\lim_{x \to x_0} \frac{|f(x)|}{|g(x)|^{1-\varepsilon}} = 0$ for all $\varepsilon > 0$.

From a results point of view, this local version of the Borg-Marchenko uniqueness theorem is our most significant new result, but a major thrust of this paper are the new methods. Theorem 1.2 says that $q$ is determined by the asymptotics of $m(-\kappa^2)$ as $\kappa \to \infty$. We can also read off differences of the boundary condition from these asymptotics. We will also prove that

THEOREM 1.3. *Let $(q_1, b_1, h_1)$, $(q_2, b_2, h_2)$ be two potentials and suppose that*

$$(1.11) \quad b_1 = b_2 \equiv b < \infty, \quad |h_1| + |h_2| < \infty, \quad q_1(x) = q_2(x) \quad \text{on } (0, b).$$

*Then*

$$(1.12) \qquad \lim_{\kappa \to \infty} e^{2b\kappa} |m_1(-\kappa^2) - m_2(-\kappa^2)| = 4(h_1 - h_2).$$

*Conversely, if (1.12) holds for some $b < \infty$ with a limit in $(0, \infty)$, then (1.11) holds.*

*Remark.* That (1.11) implies (1.12) is not so hard to see. It is the converse that is interesting.

To understand our new approach, it is useful to recall briefly the two approaches to the inverse problem for Jacobi matrices on $\ell^2(\{0, 1, 2, \ldots, \})$ [2], [8], [18]:

$$A = \begin{pmatrix} b_0 & a_0 & 0 & 0 & \ldots \\ a_0 & b_1 & a_1 & 0 & \ldots \\ 0 & a_1 & b_2 & a_2 & \ldots \\ \ldots & \ldots & \ldots & \ldots & \ldots \end{pmatrix}$$

with $a_i > 0$. Here the $m$-function is just $(\delta_0, (A - z)^{-1} \delta_0) = m(z)$ and, more generally, $m_n(z) = (\delta_n, (A^{(n)} - z)^{-1} \delta_n)$ with $A^{(n)}$ on $\ell^2(\{n, n+1, \ldots, \})$ obtained by truncating the first $n$ rows and $n$ columns of $A$. Here $\delta_n$ is the Kronecker vector, that is, the vector with 1 in slot $n$ and 0 in other slots. The fundamental theorem in this case is that $m(z) \equiv m_0(z)$ determines the $b_n$'s and $a_n$'s.

$m_n(z)$ obeys an analog of the Riccati equation (1.8):

$$\text{(1.13)} \qquad a_n^2 m_{n+1}(z) = b_n - z - \frac{1}{m_n(z)}.$$

One solution of the inverse problem is to turn (1.13) around to see that

$$\text{(1.14)} \qquad m_n(z)^{-1} = -z + b_n - a_n^2 m_{n+1}(z)$$

which, first of all, implies that as $z \to \infty$, $m_n(z) = -z^{-1} + O(z^{-2})$; so (1.14) implies

$$\text{(1.15)} \qquad m_n(z)^{-1} = -z + b_n + a_n^2 z^{-1} + O(z^{-2}).$$

Thus, (1.15) for $n = 0$ yields $b_0$ and $a_0^2$ and so $m_1(z)$ by (1.13), and then an obvious induction yields successive $b_k$, $a_k^2$, and $m_{k+1}(z)$.

A second solution involves orthogonal polynomials. Let $P_n(z)$ be the eigensolutions of the formal $(A - z)P_n = 0$ with boundary conditions $P_{-1}(z) = 0$, $P_0(z) = 1$. Explicitly,

$$\text{(1.16)} \qquad P_{n+1}(z) = a_n^{-1}[(z - b_n)P_n(z)] - a_{n-1}P_{n-1}.$$

Let $d\rho(x)$ be the spectral measure for $A$ and vector $\delta_0$ so that

$$\text{(1.17)} \qquad m(z) = \int \frac{d\rho(x)}{x - z}.$$

Then one can show that

$$\text{(1.18)} \qquad \int P_n(x)P_m(x)\,d\mu(x) = \delta_{nm}, \qquad n, m = 0, 1, \dots.$$

Thus, $P_n(z)$ is a polynomial of degree $n$ with positive leading coefficients determined by (1.18). These orthonormal polynomials are determined via Gram-Schmidt from $\rho$ and by (1.17) from $m$. Once one has the $P_n$, one can determine the $a$'s and $b$'s from the equation (1.16).

Of course, these approaches via the Riccati equation and orthogonal polynomials are not completely disjoint. The Riccati solution gives the $a_n$'s and $b_n$'s as continued fractions. The connection between continued fractions and orthogonal polynomials goes back a hundred years to Stieltjes' work on the moment problem [18].

The Gel'fand-Levitan-Marchenko [7], [11], [12], [13] approach to the continuum case is a direct analog of this orthogonal polynomial case. One looks at solutions $U(x, k)$ of

$$\text{(1.19)} \qquad -U'' + q(x)U = k^2 U(x)$$

obeying $U(0) = 1$, $U'(0) = ik$, and proves that they obey a representation

$$\text{(1.20)} \qquad U(x, k) = e^{ikx} + \int_{-x}^{x} K(x, y)e^{iky}\,dy,$$

the analog of $P_n(z) = cz^n +$ lower order. One defines $s(x,k) = (2ik)^{-1}[U(x,k) - U(x,-k)]$ which obeys (1.19) with $s(0) = 0$, $s'(0) = 1$.

The spectral measure $d\rho$ associated to $m(z)$ by

$$d\rho(\lambda) = \lim_{\varepsilon\downarrow 0}[(2\pi)^{-1}\operatorname{Im} m(\lambda + i\varepsilon)\,d\lambda]$$

obeys

(1.21) $$\int s(x,k)s(y,k)\,d\rho(k^2) = \delta(x-y),$$

at least formally. (1.20) and (1.21) yield an integral equation for $K$ depending only on $d\rho$ and then once one has $K$, one can find $U$ and so $q$ via (1.19) (or via another relation between $K$ and $q$).

Our goal in this paper is to present a new approach to the continuum case, that is, an analog of the Riccati equation approach to the discrete inverse problem. The simple idea for this is attractive but has a difficulty to overcome. $m(z,x)$ determines $q(x)$ at least if $q$ is continuous by the known asymptotics ([4]):

(1.22) $$m(-\kappa^2, x) = -\kappa - \frac{q(x)}{2\kappa} + o(\kappa^{-1}).$$

We can therefore think of (1.8) with $q$ defined by (1.22) as an evolution equation for $m$. The idea is that using a suitable underlying space and uniqueness theorem for solutions of differential equations, (1.8) should uniquely determine $m$ for all positive $x$, and so $q(x)$ by (1.22).

To understand the difficulty, consider a potential $q(x)$ on the whole real line. There are then functions $u_\pm(x,z)$ defined for $z \in \mathbb{C}\setminus[\beta,\infty)$ which are $L^2$ at $\pm\infty$ and two $m$-functions $m_\pm(z,x) = \frac{u'_\pm(x,z)}{u_\pm(x,z)}$. Both obey (1.8), yet $m_+(0,z)$ determines and is determined by $q$ on $(0,\infty)$ while $m_-(0,z)$ has the same relation to $q$ on $(-\infty,0)$. Put differently, $m_+(0,z)$ determines $m_+(x,z)$ for $x > 0$ but not at all for $x < 0$. $m_-$ is the reverse. So uniqueness for (1.8) is one-sided and either side is possible! That this does not make the scheme hopeless is connected with the fact that $m_-$ does not obey (1.22); rather

(1.23) $$m_-(-\kappa^2, x) = \kappa + \frac{q(x)}{2\kappa} + o(\kappa^{-1}).$$

We will see the one-sidedness of the solubility is intimately connected with the sign of the leading $\pm\kappa$ term in (1.22) and (1.23).

The key object in this new approach is a function $A(\alpha)$ defined for $\alpha \in (0,b)$ related to $m$ by

(1.24) $$m(-\kappa^2) = -\kappa - \int_0^a A(\alpha)e^{-2\alpha\kappa}\,d\alpha + \tilde{O}(e^{-2a\kappa})$$

as $\kappa \to \infty$. We have written $A(\alpha)$ as a function of a single variable but we will allow similar dependence on other variables. Since $m(-\kappa^2, x)$ is also an $m$-function, (1.24) has an analog with a function $A(\alpha, x)$. We will also sometimes consider the $q$-dependence explicitly, using $A(\alpha, x; q)$ or for $\lambda$ real and $q$ fixed $A(\alpha, x; \lambda) \equiv A(\alpha, x; \lambda q)$. If we are interested in $q$-dependence but not $x$, we will sometimes use $A(\alpha; \lambda)$. The semicolon and context distinguish between $A(\alpha, x)$ and $A(\alpha; \lambda)$.

By uniqueness of inverse Laplace transforms (see Theorem A.2.2 in Appendix 2), (1.24) and $m$ near $-\infty$ uniquely determine $A(\alpha)$.

Not only will (1.24) hold but, in a sense, $A(\alpha)$ is close to $q(\alpha)$. Explicitly, in Section 3 we will prove that

THEOREM 1.4. *Let $m$ be the $m$-function of the potential $q$. Then there is a function $A(\alpha) \in L^1(0,b)$ if $b < \infty$ and $A(\alpha) \in L^1(0,a)$ for all $a < \infty$ if $b = \infty$ so that (1.24) holds for any $a \leq b$ with $a < \infty$. $A(\alpha)$ only depends on $q(y)$ for $y \in [0, \alpha]$. Moreover, $A(\alpha) = q(\alpha) + E(\alpha)$ where $E(\alpha)$ is continuous and obeys*

$$(1.25) \qquad |E(\alpha)| \leq \left(\int_0^\alpha |q(y)|\, dy\right)^2 \exp\left(\alpha \int_0^\alpha |q(y)|\, dy\right).$$

Restoring the $x$-dependence, we see that $A(\alpha, x) = q(\alpha + x) + E(\alpha, x)$ where

$$\lim_{\alpha \downarrow 0} \sup_{0 \leq x \leq a} |E(\alpha, x)| = 0$$

for any $a > 0$; so

$$(1.26) \qquad \lim_{\alpha \downarrow 0} A(\alpha, x) = q(x),$$

where this holds in general in $L^1$ sense. If $q$ is continuous, (1.26) holds pointwise. In general, (1.26) will hold at any point of right Lebesgue continuity of $q$.

Because $E$ is continuous, $A$ determines any discontinuities or singularities of $q$. More is true. If $q$ is $C^k$, then $E$ is $C^{k+2}$ in $\alpha$, and so $A$ determines $k^{\text{th}}$ order kinks in $q$. Much more is true. In Section 7, we will prove

THEOREM 1.5. *$q$ on $[0, a]$ is only a function of $A$ on $[0, a]$. Explicitly, if $q_1, q_2$ are two potentials, let $A_1, A_2$ be their $A$-functions. If $a < b_1$, $a < b_2$, and $A_1(\alpha) = A_2(\alpha)$ for $\alpha \in [0, a]$, then $q_1(x) = q_2(x)$ for $x \in [0, a]$.*

Theorems 1.4 and 1.5 immediately imply Theorem 1.2. For by Theorem A.2.2, (1.10) is equivalent to $A_1(\alpha) = A_2(\alpha)$ for $\alpha \in [0, a]$. Theorems 1.4 and 1.5 say this holds if and only if $q_1(x) = q_2(x)$ for $x \in [0, a]$.

As noted, the singularities of $q$ come from singularities of $A$. A boundary condition is a kind of singularity, so one might hope that boundary conditions correspond to very singular $A$. In essence, we will see that this is the case — there are delta-function and delta-prime singularities at $\alpha = b$. Explicitly, in Section 5, we will prove that

THEOREM 1.6.  *Let $m$ be the m-function for a potential $q$ with $b < \infty$. Then for $a < 2b$,*

$$(1.27) \quad m(-\kappa^2) = -\kappa - \int_0^a A(\alpha) e^{-2\alpha\kappa} \, d\alpha - A_1 \kappa e^{-2\kappa b} - B_1 e^{-2\kappa b} + \tilde{O}(e^{-2a\kappa}),$$

*where*

(a) *If $h = \infty$, then $A_1 = 2$,   $B_1 = -2 \int_0^b q(y) \, dy$*

(b) *If $|h| < \infty$, then $A_1 = -2$,   $B_1 = 2[2h + \int_0^b q(y) \, dy]$.*

As we will see in Section 5, this implies Theorem 1.3.

The reconstruction theorem, Theorem 1.5, depends on the differential equation that $A(\alpha, x)$ obeys. Remarkably, $q$ drops out of the translation of (1.8) to the equation for $A$:

$$(1.28) \quad \frac{\partial A(\alpha, x)}{\partial x} = \frac{\partial A(\alpha, x)}{\partial \alpha} + \int_0^\alpha A(\beta, x) A(\alpha - \beta, x) \, d\beta.$$

If $q$ is $C^1$, the equation holds in the classical sense. For general $q$, it holds in a variety of weaker senses. Either way, $A(\alpha, 0)$ for $\alpha \in [0, a]$ determines $A(\alpha, x)$ for all $x, \alpha$ with $\alpha > 0$ and $0 < x + \alpha < a$. (1.26) then determines $q(x)$ for $x \in [0, a)$. That is the essence of where uniqueness comes from.

Here is a summary of the rest of this paper. In Section 2, we start the proof of Theorem 1.4 by considering $b = \infty$ and $q \in L^1(0, \infty)$. In that case, we prove a version of (1.24) with no error; namely, $A(\alpha)$ is defined on $(0, \infty)$ obeying

$$|A(\alpha) - q(\alpha)| \leq \|q\|_1^2 \exp(\alpha \|q\|_1)$$

and if $\kappa > \frac{1}{2}\|q\|_1$, then

$$(1.29) \quad m(-\kappa^2) = -\kappa - \int_0^\infty A(\alpha) e^{-2\alpha\kappa} \, d\alpha.$$

In Section 3, we use this and localization estimates from Appendix 1 to prove Theorem 1.4 in general. Section 4 is an aside to study implications of (1.24) for asymptotic expansions. In particular, we will see that

$$(1.30) \quad m(-\kappa^2) = -\kappa - \int_0^a q(\alpha) e^{-2\alpha\kappa} \, d\alpha + o(\kappa^{-1}),$$

which is essentially a result of Atkinson [1]. In Section 5, we turn to proofs of Theorems 1.6 and 1.3. Indeed, we will prove an analog of (1.27) for any $a < \infty$. If $a < nb$, then there are terms $\sum_{m=1}^n (A_m \kappa e^{-2m\kappa b} + B_m e^{-2m\kappa b})$ with explicit $A_m$ and $B_m$.

In Section 6, we prove (1.28), the evolution equation for $A$. In Section 7, we prove the fundamental uniqueness result, Theorem 1.5. Section 8 includes various comments including the relation to the Gel′fand-Levitan approach and a discussion of further questions raised by this approach.

I thank P. Deift, I. Gel′fand, R. Killip, and especially F. Gesztesy, for useful comments, and M. Ben-Artzi for the hospitality of Hebrew University where part of this work was done.

## 2. Existence of $A$: The $L^1$ case

In this section, we prove that when $q \in L^1$, then (1.29), which is a strong version of (1.24), holds. Indeed, we will prove

THEOREM 2.1.  *Let $q \in L^1(0, \infty)$. Then there exists a function $A(\alpha)$ on $(0, \infty)$ with $A - q$ continuous, obeying*

$$(2.1) \qquad |A(\alpha) - q(\alpha)| \leq Q(\alpha)^2 \exp(\alpha Q(\alpha)),$$

*where*

$$(2.2) \qquad Q(\alpha) \equiv \int_0^\alpha |q(y)|\, dy;$$

*thus if $\kappa > \frac{1}{2}\|q\|_1$, then*

$$(2.3) \qquad m(-\kappa^2) = -\kappa - \int_0^\infty A(\alpha) e^{-2\alpha\kappa}\, d\alpha.$$

*Moreover, if $q, \tilde{q}$ are both in $L^1$, then*

$$(2.4) \quad |A(\alpha; q) - A(\alpha; \tilde{q})| \leq \|q - \tilde{q}\|_1 [Q(\alpha) + \tilde{Q}(\alpha)] \exp(\alpha[Q(\alpha) + \tilde{Q}(\alpha)]).$$

We begin the proof with several remarks. First, since $m(-\kappa^2)$ is analytic in $\mathbb{C}\setminus[\beta, \infty)$, we need only prove (2.3) for all sufficiently large $\kappa$. Second, since $m(-\kappa^2; q_n) \to m(-\kappa^2; q)$ as $n \to \infty$ if $\|q_n - q\|_1 \to 0$, we can use (2.4) to see that it suffices to prove the theorem if $q$ is a continuous function of compact support, which we do henceforth. So suppose $q$ is continuous and supported in $[0, B]$.

We will prove the following:

LEMMA 2.2. *Let $q$ be a continuous function supported on $[0, B]$. For $\lambda \in \mathbb{R}$, let $m(z; \lambda)$ be the m-function for $\lambda q$. Then for any $z \in \mathbb{C}$ with $\operatorname{dist}(z, [0, \infty)) > \lambda \|q\|_\infty$,*

$$(2.5\text{a}) \qquad m(z; \lambda) = -\kappa - \sum_{n=1}^{\infty} M_n(z; q) \lambda^n,$$

*where for $\kappa > 0$,*

$$(2.5\text{b}) \qquad M_n(-\kappa^2; q) = \int_0^{nB} A_n(\alpha) e^{-2\kappa\alpha} \, d\alpha,$$

*where*

$$(2.6) \qquad A_1(\alpha) = q(\alpha)$$

*and for $n \geq 2$, $A_n(\alpha)$ is a continuous function obeying*

$$(2.7) \qquad |A_n(\alpha)| \leq Q(\alpha)^n \frac{\alpha^{n-2}}{(n-2)!}.$$

*Moreover, if $\tilde{q}$ is a second such potential and $n \geq 2$,*

$$(2.8) \quad |A_n(\alpha; q) - A_n(\alpha; \tilde{q})| \leq (Q(\alpha) + \tilde{Q}(\alpha))^{n-1} \left[ \int_0^\alpha |q(y) - \tilde{q}(y)| \, dy \right] \frac{\alpha^{n-2}}{(n-2)!}.$$

*Proof of Theorem 2.1 given Lemma 2.2.* By (2.7),

$$\int_0^\infty \sum_{n=2}^\infty |A_n(\alpha)| e^{-2\kappa\alpha} \, d\alpha < \infty$$

if $\kappa > \frac{1}{2} \|q\|_1$. Thus in (2.5a) for $\lambda = 1$, we can interchange the sum and integral to get the representation (2.3). (2.7) then implies (2.1) and (2.8) implies (2.4). □

*Proof of Lemma 2.2.* Let $H_\lambda$ be $-\frac{d^2}{dx^2} + \lambda q(x)$ on $L^2(0, \infty)$ with $u(0) = 0$ boundary conditions at 0. Then $\|(H_0 - z)^{-1}\| = \operatorname{dist}(z, [0, \infty))^{-1}$. So, in the sense of $L^2$ operators, if $\operatorname{dist}(z, [0, \infty)) > \lambda \|q\|_\infty$, the expansion

$$(2.9) \qquad (H_\lambda - z)^{-1} = \sum_{n=0}^\infty (-1)^n (H_0 - z)^{-1} [\lambda q (H_0 - z)^{-1}]^n$$

is absolutely convergent.

As is well known, $G_\lambda(x, y; z)$, the integral kernel of $(H_\lambda - z)^{-1}$, can be written down in terms of the solution $u$ which is $L^2$ at infinity, and the solution $w$ of

$$(2.10) \qquad -w'' + qw = zw$$

obeying $w(0) = 0$, $w'(0) = 1$

(2.11) $$G_\lambda(x, y; z) = w(\min(x, y)) \frac{u(\max(x, y))}{u(0)}.$$

In particular,

(2.12) $$m(z) = \lim_{\substack{x<y \\ y \downarrow 0}} \frac{\partial^2 G}{\partial x \partial y}.$$

From this and (2.9), we see that (using $\frac{\partial G_0}{\partial x}(x, y)\big|_{x=0} = e^{-\kappa y}$)

$$m(-\kappa^2; \lambda) = -\kappa - \lambda \int e^{-2\kappa y} q(y)\, dy + \lambda^2 \langle \varphi_\kappa, (H_\lambda + \kappa^2)^{-1} \varphi_\kappa \rangle,$$

where $\varphi_\kappa(y) = q(y)e^{-\kappa y}$. Since $\varphi_\kappa \in L^2$, we can use the convergent expansion (2.9) and so conclude that (2.5a) holds with (for $n \geq 2$)

(2.13) $$M_n(-\kappa^2; q) = (-1)^{n-1} \int e^{-\kappa x_1} q(x_1) G_0(x_1, x_2) q(x_2)$$
$$\ldots G_0(x_{n-1}, x_n) q(x_n) e^{-\kappa x_n}\, dx_1 \ldots dx_n.$$

Now use the following representation for $G_0$:

(2.14) $$G_0(x, y; -\kappa^2) = \frac{\sinh(\kappa \min(x, y))}{\kappa} e^{-\kappa \max(x, y)}$$
$$= \frac{1}{2} \int_{|x-y|}^{x+y} e^{-\ell \kappa}\, d\ell$$

to write

(2.15)
$$M_n(-\kappa^2; q)$$
$$= \frac{(-1)^{n-1}}{2^{n-1}} \int_{R_n} q(x_1) \ldots q(x_n) e^{-2\alpha(x_1, x_n, \ell_1, \ldots, \ell_{n-1})\kappa}\, dx_1 \ldots dx_n d\ell_1 \ldots d\ell_{n-1},$$

where $\alpha$ is shorthand for the linear function

(2.16) $$\alpha = \frac{1}{2}\left(x_1 + x_n + \sum_{j=1}^{n-1} \ell_j\right)$$

and $R_n$ is the region

$$R_n = \{(x_1, \ldots, x_n, \ell_1, \ldots, \ell_{n-1}) \in \mathbb{R}^{2n-1} \mid 0 \leq x_i \leq B \text{ for } i = 1, \ldots, n;$$
$$|x_i - x_{i+1}| \leq \ell_i \leq x_i + x_{i-1} \text{ for } i - 1, \ldots, n-1\}.$$

In the region $R_n$, notice that

$$\alpha \leq \frac{1}{2}\left(x_1 + x_n + \sum_{j=1}^{n-1}(x_j + x_{j+1})\right) = \sum_{j=1}^{n} x_j \leq nB.$$

Change variables by replacing $\ell_{n-1}$ by $\alpha$ using the linear transformation (2.16) and use $\ell_{n-1}$ for the linear function

$$\ell_{n-1}(x_1, x_n, \ell_1, \ldots, \ell_{n-2}, \alpha) = 2\alpha - x_1 - x_n - \sum_{j=1}^{n-2} \ell_j. \tag{2.17}$$

Thus, (2.5b) holds where

$$A_n(\alpha) = \frac{(-1)^{n-1}}{2^{n-2}} \int_{R_n(\alpha)} q(x_1) \ldots q(x_n) \, dx_1 \ldots dx_n d\ell_1 \ldots d\ell_{n-2}. \tag{2.18}$$

$2^{n-1}$ has become $2^{n-2}$ because of the Jacobian of the transition from $\ell_{n-1}$ to $\alpha$. $R_n(\alpha)$ is the region

$$\tag{2.19}$$
$$R_n(\alpha) = \{(x_1, \ldots, x_n, \ell_1, \ldots, \ell_{n-2}) \in \mathbb{R}^{2n-2} \mid 0 \le x_i \le B \text{ for } i = 1, \ldots, n;$$
$$|x_i - x_{i+1}| \le \ell_i \le x_i + x_{i+1} \text{ for } i = 1, \ldots, n-2;$$
$$|x_{n-1} + x_n| \le \ell_{n-1}(x_1, \ldots, x_n, \ell_1, \ldots, \ell_{n-2}, \alpha) \le x_{n-1} + x_n\}$$

with $\ell_{n-1}$ the functional given by (2.17).

We claim that

$$\tag{2.20}$$
$$R_n(\alpha) \subset \tilde{R}_n(\alpha)$$
$$= \left\{(x_1, \ldots, x_n, \ell_1, \ldots, \ell_{n-2}) \in \mathbb{R}^{2n-2} \bigg| 0 \le x_i \le \alpha; \, \ell_i \ge 0; \, \sum_{i=1}^{n-2} \ell_i \le 2\alpha\right\}.$$

Accepting (2.20) for a moment, we note by (2.18) that

$$|A_n(\alpha)| \le \frac{1}{2^{n-2}} \int_{\tilde{R}_n(\alpha)} |q(x_1)| \ldots |q(x_n)| \, dx_1 \ldots d\ell_{n-2}$$
$$= \left(\int_0^\alpha |q(x)| \, dx\right)^n \frac{\alpha^{n-2}}{(n-2)!}$$

since $\int_{\sum y_i = b; y_i \ge 0} dy_1 \ldots dy_n = \frac{b^n}{n!}$ by a simple induction. This is just (2.7).

To prove (2.8), we note that

$$|A_n(\alpha; q) - A_n(\alpha, \tilde{q})|$$
$$\le 2^{-n-2} \int_{\tilde{R}_n(\alpha)} |q(x_1) \ldots q(x_n) - \tilde{q}(x_1) \ldots \tilde{q}(x_n)| \, dx_1 \ldots d\ell_{n-2}$$
$$\le \frac{\alpha^{n-2}}{(n-2)!} \sum_{j=0}^{n-1} Q(\alpha)^j \left[\int_0^\alpha |q(y) - \tilde{q}(y)| \, dy\right] \tilde{Q}(\alpha)^{n-j-1}.$$

Since $\sum_{j=0}^m a^j b^{m-j} \le \sum_{j=0}^m \binom{m}{j} a^j b^{m-j} = (a+b)^m$, (2.8) holds.

Thus, we need only prove (2.20). Suppose $(x_1, \ldots, x_n, \ell_1, \ldots, \ell_{n-2}) \in R_n(\alpha)$. Then

$$2x_m \leq |x_1 - x_m| + |x_n - x_m| + x_1 + x_n$$
$$\leq x_1 + x_n + \sum_{j=1}^{n-1} |x_{j+1} - x_j|$$
$$\leq x_1 + x_n + \sum_{j=1}^{n-2} \ell_j + \ell_{n-1}(x_1, \ldots, x_n, \ell_1, \ldots, \ell_{n-2}; \alpha) = 2\alpha$$

so $0 \leq x_j \leq \alpha$, proving that part of the condition $(x_1, \ell_{n-2}) \subset \tilde{R}_n(\alpha)$. For the second part, note that

$$\sum_{j=1}^{n-2} \ell_j = 2\alpha - x_1 - x_n - \ell_{n-1}(x_1, \ldots, x_n, \ell_1, \ldots, \ell_{n-2}) \leq 2\alpha$$

since $x_1$, $x_n$, and $\ell_{n-2}$ are nonnegative on $R_n(\alpha)$. □

We want to say more about the smoothness of the functions $A_n(\alpha)$ and $A_n(\alpha, x)$ defined for $x \geq 0$ and $n \geq 2$ by

(2.21)
$$A_n(\alpha, x) = \frac{(-1)^{n-1}}{2^{n-2}} \int_{R_n(\alpha)} q(x + x_1) \ldots q(x + x_n) \, dx_1 \ldots dx_n d\ell_1 \ldots d\ell_{n-2}$$

so that $A(\alpha, x) = \sum_{n=0}^{\infty} A_n(\alpha, x)$ is the $A$-function associated to $m(-\kappa^2, x)$. We begin with $\alpha$ smoothness for fixed $x$.

PROPOSITION 2.3.  $A_n(\alpha, x)$ is a $C^{n-2}$-function in $\alpha$ and obeys for $n \geq 3$

(2.22) $$\left| \frac{d^j A_n(\alpha)}{d\alpha^j} \right| \leq \frac{1}{(n-2-j)!} \alpha^{n-2-j} Q(\alpha)^n; \qquad j = 1, \ldots, n-2.$$

*Proof.* Write

$$A_n(\alpha) = \frac{(-1)^{n-1}}{2^{n-1}} \int_{R_n} q(x_1) \ldots q(x_n) \delta\left(2\alpha - x_1 - x_n - \sum_{m=1}^{n-1} \ell_i\right) dx_1 \ldots dx_n d\ell_j \ldots d\ell_{n-1}.$$

Thus, formally,

(2.23)
$$\frac{d^j A_n(\alpha)}{d\alpha^j} = \frac{(-1)^{n-1} 2^j}{2^{n-2}} \int_{R_n} q(x_1) \ldots q(x_n) \delta^{(j)}\left(2\alpha - x_1 - x_n - \sum_{m=1}^{n-1} \ell_i\right) dx_1 \ldots d\ell_{n-1}.$$

Since $j + 1 \leq n - 1$, we can successively integrate out $\ell_{n-1}, \ell_{n-2}, \ldots, \ell_{n-j-1}$ using

$$(2.24) \qquad \int_a^b \delta^j(c - \ell)\, d\ell = \delta^{j-1}(c - a) - \delta^{j-1}(c - b)$$

and

$$(2.25) \qquad \int_a^b \delta(c - \ell)\, d\ell = \chi_{(a,b)}(c).$$

Then we estimate each of the resulting $2^j$ terms as in the previous lemma, getting

$$\left|\frac{d^j A_n(\alpha)}{d\alpha^j}\right| \leq \frac{2^j}{2^{n-2}}\, Q(\alpha)^n \frac{(2\alpha)^{n-j-2}}{(n-j-2)!}$$

which is (2.22).

(2.24), (2.25), while formal, are a way of bookkeeping for legitimate movement of hyperplanes. In (2.25), there is a singularity at $c = a$ and $c = b$, but since we are integrating in further variables, these are irrelevant. $\square$

PROPOSITION 2.4. *If $q$ is $C^m$, then $A_n(\alpha)$ is $C^{m+(2n-2)}$.*

*Proof.* Write $R_n$ as $n!$ terms with orderings $x_{\pi(1)} < \cdots < x_{\pi(n)}$. For $j_0 = 2n - 2$, we integrate out all $2n - 1$, $\ell$ and $x$ variables. We get a formula for $\frac{d^{j_0} A_n(\alpha)}{d\alpha^{j_0}}$ as a sum of products of $q$'s evaluated at rational multiples of $\alpha$. We can then take $m$ additional derivatives. $\square$

THEOREM 2.5. *If $q$ is $C^m$ and in $L^1(0, \infty)$, then $A(\alpha)$ is $C^m$ and $A(\alpha) - q(\alpha)$ is $C^{m+2}$.*

*Proof.* By (2.2), we can sum the terms in the series for $\frac{d^j A}{d\alpha^j}$ and $\frac{d^j(A-q)}{d\alpha^j}$ for $j = 0, 1, \ldots, m$ and $j = 0, 1, \ldots, m-2$, respectively. With this bound and the fundamental theorem of calculus, one can prove the stated regularity. $\square$

Now we can turn to $x$-dependence.

LEMMA 2.6. *If $q$ is $C^k$ and of compact support, then $A_n(\alpha, x)$ for $\alpha$ fixed is $C^k$ in $x$, and for $n \geq 2$, $j = 1, \ldots, k$,*

$$(2.26) \qquad \left|\frac{d^j A_n(\alpha, x)}{dx^j}\right| \leq Q(\alpha)^{\max(0, n-j)} [P_j(\alpha)]^{\min(j,n)} \frac{\alpha^{n-2}}{(n-2)!},$$

*where*

$$P_j(\alpha) = \int_0^\alpha \sum_{m=0}^j \left|\frac{d^m q}{dx^m}(y)\right| dy.$$

*Proof.* In (2.21), we can take derivatives with respect to $x$. We get a sum of terms with derivatives on each $q$, and using values on these terms and the argument in the proof of Lemma 2.2, we obtain (2.26). □

THEOREM 2.7. *If $q$ is $C^k$ and of compact support, then $A(\alpha, x)$ for $\alpha$ fixed is $C^k$ in $x$ and*

$$\frac{d^j m}{dx^j}(-\kappa^2, x) = -\int_0^\infty \frac{\partial^j A}{\partial x^j}(\alpha, x) e^{-2\alpha\kappa}\, d\alpha$$

*for $\kappa$ large and $j = 1, 2, \ldots, k$.*

*Proof.* This follows from the estimates in Lemma 2.6 and Theorem 2.1. □

## 3. Existence of $A$: General case

By combining Theorem 2.1 and Theorem A.1.1, we immediately have

THEOREM 3.1. *Let $b < \infty$, $q \in L^1(0, b)$, and $h \in \mathbb{R} \cup \{\infty\}$ or else let $b = \infty$ and let $q$ obey (1.4), (1.5). Fix $a < b$. Then, there exists a function $A(\alpha)$ on $L^1(0, a)$ obeying*

(3.1) $$|A(\alpha) - q(\alpha)| \leq Q(\alpha) \exp(\alpha Q(\alpha)),$$

*where*

(3.2) $$Q(\alpha) \equiv \int_0^\alpha |q(y)|\, dy$$

*so that as $\kappa \to \infty$,*

(3.3) $$m(-\kappa^2) = -\kappa - \int_0^a A(\alpha) e^{-2\alpha\kappa}\, d\alpha + \tilde{O}(e^{-2a\kappa}).$$

*Moreover, $A(\alpha)$ on $[0, a]$ is only a function of $q$ on $[0, a]$.*

*Proof.* Let $\tilde{b} = \infty$ and $\tilde{q}(x) = q(x)$ for $x \in [0, a]$ and $\tilde{q}(x) = 0$ for $x > a$. By Theorem A.1.1, $m - \tilde{m} = \tilde{O}(e^{-2a\kappa})$, and by Theorem 2.1, $\tilde{m}$ has a representation of the form (3.3). □

## 4. Asymptotic formula

While our interest in the representation (1.24) is primarily for inverse theory and, in a sense, it provides an extremely complete form of asymptotics, the formula is also useful to recover and extend results of others on more conventional asymptotics.

In this section, we will explain this theme. We begin with a result related to Atkinson [1] (who extended Everitt [5]).

THEOREM 4.1. *For any $q$ (obeying (1.2)–(1.5)), we have that*

$$(4.1) \qquad m(-\kappa^2) = -\kappa - \int_0^b q(x) e^{-2x\kappa}\, dx + o(\kappa^{-1}).$$

*Remarks.* 1. Atkinson's "$m$" is the negative inverse of our $m$ and he uses $k = i\kappa$, and so his formula reads ((4.3) in [1])

$$m_{\text{Atk}}(k^2) = ik^{-1} + k^{-2} \int_0^b e^{2ikx} q(x)\, dx + o(|k|^{-3}).$$

2. Atkinson's result is stronger in that he allows cases where $q$ is not bounded below (and so he takes $|z| \to \infty$ staying away from the negative real axis also). [10] will extend (4.1) to some such situations.

3. Atkinson's method breaks down on the real $x$ axis where our estimates hold, but one could use Phragmén-Lindelöf methods and Atkinson's results to prove Theorem 4.1.

*Proof.* By Theorem 3.1, $(A - q) \to 0$ as $\alpha \downarrow 0$ so $\int_0^a e^{-2a\kappa}(A(\alpha) - q(\alpha))\, d\alpha = o(\kappa^{-1})$. Thus, (3.3) implies (4.1). □

COROLLARY 4.2.
$$m(-\kappa^2) = -\kappa + o(1).$$

*Proof.* Since $q \in L^1$, dominated convergence implies that $\int_0^b q(x) e^{-2\kappa x}\, dx = o(1)$. □

COROLLARY 4.3. *If $\lim_{x \downarrow 0} q(x) = a$ (indeed, if $\frac{1}{s} \int_0^s q(x)\, dx \to a$ as $s \downarrow 0$), then*

$$m(-\kappa^2) = -\kappa - \frac{a}{2}\kappa^{-1} + o(\kappa^{-1}).$$

COROLLARY 4.4. *If $q(x) = cx^{-\alpha} + o(x^{-\alpha})$ for $0 < \alpha < 1$, then*

$$m(-\kappa^2) = -\kappa - c[2^{a-1}\Gamma(1-\alpha)]\kappa^{\alpha-1} + o(\kappa^{\alpha-1}).$$

We can also recover the result of Danielyan and Levitan [4]:

THEOREM 4.5. *Let $q(x) \in C^n[0, \delta)$ for some $\delta > 0$. Then as $\kappa \to \infty$, for suitable $\beta_0, \ldots, \beta_n$, we have that*

$$(4.2) \qquad m(-\kappa^2) = -\kappa - \sum_{m=0}^n \beta_j \kappa^{-j-1} + O(\kappa^{-n-1}).$$

*Remarks.* 1. Our $m$ is the negative inverse of their $m$.

2. Our proof does not require that $q$ is $C^n$. It suffices that $q(x)$ has an asymptotic series $\sum_{m=0}^{n} a_m x^m + o(x^n)$ as $x \downarrow 0$.

*Proof.* By Theorems 3.1 and 2.5, $A(\alpha)$ is $C^n$ on $[0, \delta)$. It follows that $A(\alpha) = \sum_{m=0}^{n} b_j \alpha^j + o(\alpha^j)$. Since $\int_0^\delta \alpha^j e^{-2\alpha\kappa} \, d\alpha = \kappa^{-j-1} 2^{-j-1} j! + \tilde{O}(e^{-2\delta\kappa})$, we have (4.2) $\beta_j = 2^{j-1} j! b_j = 2^{j-1} \frac{\partial^j A}{\partial \alpha^j}(\alpha = 0)$. □

Later we will prove that $A$ obeys (1.28). This immediately yields a recursion formula for $\beta_j(x)$, viz.:

$$\beta_{j+1}(x) = \frac{1}{2} \frac{\partial \beta_j}{\partial x} + \frac{1}{2} \sum_{\ell=0}^{j} \beta_\ell(x) \beta_{j-\ell}(x), \qquad j \geq 0$$

$$\beta_0(x) = \frac{1}{2} q(x);$$

see also [9, §2].

## 5. Reading boundary conditions

Our goal in this section is to prove Theorem 1.6 and then Theorem 1.3. Indeed, we will prove the following stronger result:

THEOREM 5.1. *Let $m$ be the m-function for a potential $q$ with $b < \infty$. Then there exists a measurable function $A(\alpha)$ on $[0, \infty)$ which is $L^1$ on any finite interval $[0, R]$, so that for each $N = 1, 2, \ldots$ and any $a < 2Nb$,*

(5.1)
$$m(-\kappa^2) = -\kappa - \int_0^a A(\alpha) e^{-2\alpha\kappa} \, d\alpha - \sum_{j=1}^{N} A_j \kappa e^{-2\kappa bj} - \sum_{j=1}^{N} B_j e^{-2\kappa bj} + \tilde{O}(e^{-2a\kappa}),$$

*where*

(a) *If $h = \infty$, then $A_j = 2$ and $B_j = -2j \int_0^b q(y) \, dy$.*

(b) *If $|h| < \infty$, then $A_j = 2(-1)^j$ and $B_j = 2(-1)^{j+1} j [2h + \int_0^b q(y) \, dy]$.*

*Remarks.* 1. The combination $2h + \int_0^b q(y) \, dy$ is natural when $|h| < \infty$. It also enters into the formula for eigenvalue asymptotics [11], [13].

2. One can think of (5.1) as saying that

$$m(-\kappa^2) = -\kappa - \int_0^a \tilde{A}(\alpha) e^{-2\alpha\kappa} \, d\alpha + \tilde{O}(e^{-2a\kappa})$$

for any $a$ where now $\tilde{A}$ is only a distribution of the form $\tilde{A}(\alpha) = A(\alpha) + \frac{1}{2} \sum_{j=1}^{\infty} A_j \delta'(\alpha - jb) + \sum_{j=1}^{\infty} B_j \delta(\alpha - jb)$ where $\delta'$ is the derivative of a delta function.

3. As a consistency check on our arithmetic, we note that if $q(y) \to q(y)+c$ and $\kappa^2 \to \kappa^2 - c$ for some $c$, then $m(-\kappa^2)$ should not change. $\kappa^2 \to \kappa^2 - c$ means $\kappa \to \kappa - \frac{c}{2\kappa}$ and so $\kappa e^{-2\kappa bj} \to \kappa e^{-2\kappa bj} + cbj e^{-2\kappa bj} + O(\kappa^{-1})$ terms. That means that under $q \to q + c$, we must have that $B_j \to B_j - cbj A_j$, which is the case.

*Proof.* Consider first the free Green's function for $-\frac{d^2}{dx^2}$ with Dirichlet boundary conditions at $0$ and $h$-boundary condition at $b$. It has the form

$$(5.2) \qquad G_0(x,y) = \frac{\sinh(\kappa x)\, u_+(y)}{\kappa\, u_+(0)}, \qquad x < y$$

where $u_+(y; \kappa, h)$ obeys $-u'' = -\kappa^2 u$ with boundary condition

$$(5.3) \qquad u'(b) + hu(b) = 0.$$

Write

$$(5.4) \qquad u_+(y) = e^{-\kappa y} + \alpha e^{-\kappa(2b-y)}$$

for $\alpha \equiv \alpha(h, \kappa)$. Plugging (5.4) into (5.3), one finds that

$$(5.5) \qquad \alpha = \begin{cases} -1, & h = \infty \\ \frac{1-h/\kappa}{1+h/\kappa} = 1 - \frac{2h}{\kappa} + O(\kappa^{-2}), & |h| < \infty. \end{cases}$$

Now one just follows the arguments of Section 2 using (5.2) in place of (2.14). All terms of order 2 or more in $\lambda^2$ contribute to locally $L^1$ pieces of $\tilde{A}(\alpha)$. The exceptions come from the order 0 and order 1 terms. The order 0 term is

$$\lim_{x < y \to 0} \frac{\partial^2 G_0(x,y)}{\partial x \partial y} = \frac{u'_+(0)}{u_+(0)} = -\kappa \left[ \frac{1 - \alpha e^{-2b\kappa}}{1 + \alpha e^{-2b\kappa}} \right] \equiv Q.$$

Now $\frac{1-z}{1+z} = 1 + 2 \sum_{n=1}^{\infty} (-1)^n z^n$, so

$$(5.6)$$
$$Q = -\kappa - 2\kappa \sum_{n=1}^{\infty} (-1)^n \alpha^n e^{-2b\kappa n}$$
$$= \begin{cases} -\kappa - 2\kappa \sum_{n=1}^{\infty} e^{-2b\kappa n} \\ -\kappa - 2\kappa \sum_{n=1}^{\infty} (-1)^n e^{-2b\kappa n} - 4 \sum_{n=1}^{\infty} (-1)^{n+1} nh e^{-2b\kappa n} + \text{regular}, \end{cases}$$

where "regular" means a term which is a Laplace transform of a locally $L^1$ function. We used (by (5.5)) that if $h$ is finite, then

$$\alpha^n = 1 - \frac{2nh}{\kappa} + O(\kappa^{-2}),$$

where $\kappa O(\kappa^{-2})$ in this context is regular.

The first-order term is

$$P \equiv -\int_0^b q(y) \left[ \frac{u_+(y)}{u_+(0)} \right]^2 dy.$$

Now $(\frac{u_+(y)}{u_+(0)})^2 = (1 + \alpha e^{-2b\kappa})^{-2}[e^{-\kappa y} + \alpha e^{-\kappa(2b-y)}]^2$. In expanding the last square, $e^{-2\kappa y}$ and $e^{-2\kappa(2b-y)}$ yield regular terms but the cross term is not regular; that is,

$$P = -\left[\int_0^b q(y)\,dy\right] 2\alpha e^{-2\kappa b}(1 + \alpha e^{-2\kappa b})^{-2} + \text{regular}.$$

Now

$$z(1+z)^{-2} = -z\frac{d}{dz}(1+z)^{-1} = -z\frac{d}{dz}\left(\sum_{n=0}^{\infty}(-1)^n z^n\right) = \sum_{n=1}^{\infty}(-1)^{n+1} n z^n$$

and so using $\alpha^n = (-1)^n$ if $h = \infty$ and $\alpha^n = 1 + O(\kappa^{-1})$ if $h < \infty$, we see that

(5.7) $$P = \begin{cases} 2\sum_{n=1}^{\infty} n e^{-2n\kappa b}[\int_0^b q(y)\,dy] + \text{regular}, & h = \infty \\ 2\sum_{n=1}^{\infty}(-1)^n n e^{-2n\kappa b}[\int_0^b q(y)\,dy] + \text{regular}, & |h| < \infty. \end{cases}$$

Combining (5.6) and (5.7), we see that (with $I = \int_0^b q(y)\,dy$),
(5.8)
$$P + Q = \begin{cases} -\kappa - 2\kappa\sum_{n=1}^{\infty} e^{-2b\kappa n} + 2\sum_{n=1}^{\infty} nI e^{-2b\kappa n} + \text{regular} \\ -\kappa - 2\kappa\sum_{n=1}^{\infty}(-1)^n e^{-2b\kappa n} + 2\sum_{n=1}^{\infty}(-1)^n n[I + 2h]e^{-2b\kappa h} + \text{regular}. \end{cases}$$

This is precisely what conclusion (a), (b) of Theorem 5.1 asserts. □

*Proof of Theorem* 1.3. The direct assertion follows from Theorem 5.1 and the fact that $A$ on $[0, b]$ is only a function of $q$ there. We consider the converse part. By Theorems 5.1 and 3.1, for each $q_j$ we have for any $a < \infty$,

$$m_j(-\kappa^2) = -\kappa - \int_0^a \tilde{A}_j(\alpha)e^{-2\kappa\alpha}\,d\alpha + \tilde{O}(e^{-2\kappa a}),$$

where $\tilde{A}(\alpha)$ is an $L^1(0, a)$ function plus a possible finite sum of $\delta$ and $\delta'$ terms. Take $a = 2b$. (1.12) and the fundamental expansion on uniqueness of inverse Laplace transforms (see Theorem A.2.2) imply that $(A_1 - A_2)(\alpha)$ is supported on $[b, 2b]$. If $b_1, b_2 > b$, then the limit (1.12) is zero, so $h_1 \neq h_2$ implies either $b_1$ or $b_2$ is $b$. If only one is $b$, then the difference has a $\delta'$ term and the limit in (1.12) is infinite. Therefore, $b_1 = b_2 = b$.

Since $A_1 = A_2$ on $[0, b]$, Theorem 1.2 implies that $q_1(x) = q_2(x)$ on $[0, b]$. If both $h_1$ and $h_2$ are infinite, then the limit is zero. If only one is infinite, then there is a $\delta'$ term and the limit is infinite. Thus, a limit on $(0, \infty)$ implies $h_1$ and $h_2$ are both finite and so, by Theorem 5.1, the limit is $4(h_1 - h_2)$ as claimed. □

### 6. The $A$-equation

In this section, we will prove equation (1.28). We begin with the case where $q$ is $C^1$. In general, given $q$ (i.e., $q$, $b$, and $h$ if $b < \infty$), we can define $m(z,x) = \frac{u'_+(x,z)}{u_+(x,z)}$ for $x \in [0,b)$ and $z \in \mathbb{C}\backslash[\beta, \infty)$ for suitable $\beta \in \mathbb{R}$. By Theorem 3.1, there is a function $A(\alpha, x)$ defined for $(\alpha, x) \in \{(\alpha, x) \in \mathbb{R}^2 \mid 0 \leq x < b; 0 < \alpha < b - x\} \equiv S$ so that for any $a < b - x$,

$$(6.1) \qquad m(-\kappa^2, x) = -\kappa - \int_0^a A(\alpha, x) e^{-2\alpha\kappa} \, d\alpha + \tilde{O}(e^{-2a\kappa}).$$

Moreover, $m$ obeys the Riccati equation (1.8), and by (3.1) if we define $g_\alpha(x)$ on $[0, b]$ by

$$g_\alpha(x) = A(\alpha, x) \quad \text{if } x < b - \alpha$$
$$= 0 \quad \text{if } b - \alpha \leq x < b,$$

then

$$(6.2) \qquad \lim_{\alpha \downarrow 0} g_\alpha(x) = q(x)$$

in $L^1(0, a)$ for any $a < b$.

In (6.2), there is a potential difficulty in that $A(\alpha, x)$ is *a priori* only defined for almost every $\alpha$ for each $x$, so that $g_\alpha(x)$ is not well-defined for all $\alpha$. One can finesse this difficulty by interpreting (6.2) in essential sense (i.e., for all $a < b$ and $\varepsilon > 0$, there is a $\Lambda$ so that for almost every $\alpha$ with $0 < \alpha < \Lambda$, we have $\int_0^a |g_\alpha(x) - q(x)| \, dx < \varepsilon$). Alternatively, one can pick a concrete realization of $q$ and then use the fact that $A - q$ is continuous to define $A(x, \alpha) - q(x + \alpha)$ for all $x, \alpha$ and then (6.2) holds in traditional sense. Indeed, if $q$ is continuous, it holds pointwise.

THEOREM 6.1.  *If $q$ is $C^1$, then $A$ is jointly $C^1$ on $S$ and obeys*

$$(6.3) \qquad \frac{\partial A}{\partial x} = \frac{\partial A}{\partial \alpha} + \int_0^\alpha A(\beta, x) A(\alpha - \beta, x) \, d\beta.$$

*Proof.* That $A$ is jointly $C^1$ when $q$ is $C^1$ of compact support follows from the arguments in Section 2 (and then the fact that $A$ on $[0, a)$ is only a function of $q$ on $[0, a)$ lets us extend this to all $C^1$ $q$'s). Moreover, by Theorem 2.7,

$$(6.4) \qquad \frac{\partial m}{\partial x}(-\kappa^2, x) = -\int_0^a \frac{\partial A}{\partial x}(\alpha, x) e^{-2\alpha\kappa} \, d\alpha + \tilde{O}(e^{-2a\kappa})$$

for all $a < b - x$. Now in (6.1), square $m$ to see that
(6.5)
$$m(x, -\kappa^2)^2 = \kappa^2 + \int_0^a B(\alpha, x) e^{-2\alpha\kappa} \, d\alpha + 2\int_0^a A(\alpha, x) \kappa e^{-2\alpha\kappa} \, d\alpha + \tilde{O}(e^{-2a\kappa}),$$

where $B(\alpha, x) = \int_0^\alpha A(\beta, x) A(\alpha - \beta, x) \, d\beta$. In the cross term in (6.5), write $2\kappa e^{-2\alpha\kappa} = -\frac{d}{d\alpha}(e^{-2\alpha\kappa})$ and integrate by parts

$$2 \int_0^a A(\alpha, x) \kappa e^{-2\alpha\kappa} \, d\alpha = -A(a, x) e^{-2a\kappa} + \lim_{\alpha \downarrow 0} A(\alpha, x) + \int_0^a \frac{\partial A}{\partial \alpha}(\alpha, x) e^{-2\alpha\kappa} \, d\alpha.$$

By (6.2), $\lim_{\alpha \downarrow 0} A(\alpha, x) = q(x)$ so (6.5) becomes

$$(6.6) \qquad -m^2 + \kappa^2 + q = \int_0^a \left(\frac{\partial A}{\partial \alpha} + B\right) e^{-2\alpha\kappa} \, d\alpha + \tilde{O}(e^{-2a\kappa}).$$

The Riccati equation (1.8), (6.4), (6.6), and the uniqueness of inverse Laplace transforms (Theorem A.2.2) then imply that (6.3) holds pointwise.  □

There are various senses in which (6.3) holds for general $q$. We will state three. All follow directly from the regularity results in Section 2, the continuity expressed by (3.4), and Theorem 6.1.

THEOREM 6.2. *For general $q$, (6.3) holds in distributional sense.*

THEOREM 6.3. *For general $q$, define $C(\gamma, x)$ on $\{(\gamma, x) \in \mathbb{R}^2 \mid x < \gamma < b)\}$ by*

$$C(\gamma, x) = A(\gamma - x, x).$$

*Then, if $x_1 < x_2 < \gamma$, we have that for all $(\gamma, x)$,*

$$(6.7) \qquad C(\gamma, x_2) = C(\gamma, x_1) + \int_{x_1}^{x_2} dy \left[\int_y^\gamma C(\lambda, y) C(\gamma - \lambda + y, y) \, d\lambda\right].$$

THEOREM 6.4. *If $q$ is continuous, then $F(\alpha, x) \equiv A(\alpha, x) - q(\alpha + x)$ is jointly $C^1$ and obeys*

$$\frac{\partial F}{\partial x} = \frac{\partial F}{\partial \alpha} + \int_0^\alpha A(\beta, x) A(\alpha - \beta, x) \, d\beta.$$

## 7. The uniqueness theorem

In this section, we will prove Theorem 1.5 and therefore, as already noted in the introduction, Theorem 1.2. Explicitly,

THEOREM 7.1. *Let $q_1$ and $q_2$ be two potentials and let $a < \min(b_1, b_2)$. Suppose $A_1(\alpha, 0) = A_2(\alpha, 0)$ for $\alpha \in [0, a]$. Then $q_1 = q_2$ for a.e. for $x$ in $[0, a]$.*

*Proof.* We will use (6.7) and an elementary Gronwall's equality to conclude that $A_1(\alpha, x) = A_2(\alpha, x)$ on $S = \{(x, \alpha) \in \mathbb{R}^2 \mid x + \alpha < a\}$, and then conclude that $q_1 = q_2$ on $[0, a]$ by (6.2). Pick an explicit realization of $q_1$ and $q_2$ and then since $A_j(\alpha, x) - q_j(\alpha + x)$ is continuous, an explicit realization of $A_j(\alpha, x)$ in which

$$g(x) = \int_0^{a-x} |A_1(\alpha, x) - A_2(\alpha, x)| \, d\alpha$$

is continuous. Moreover, in this realization,

$$D = \sup_{0 \leq x < a} \int_0^{a-x} \big[|A_1(\alpha, x)| + |A_2(\alpha, x)|\big] \, d\alpha < \infty$$

since the integral is also continuous. By (6.7) for $0 \leq x_1 < x_2 < a$,

(7.1) $$g(x_2) \leq g(x_1) + D \int_{x_1}^{x_2} g(y) \, dy.$$

Letting $h(x) = \sup_{0 \leq y \leq x} g(y)$, we see that (7.1) implies

$$h(x_2) \leq h(x_1) + D \int_{x_1}^{x_2} h(x_2) \, dy$$

so if $D(x_2 - x_1) < 1$ and $h(x_1) = 0$, then $h(x_2) = 0$. By hypothesis, $h(0) = 0$. So using this argument a finite number of times, $h(x) \equiv 0$ for $x \in [0, a]$, that is, $A_1 = A_2$ on $S$. □

## 8. Complements and open questions

In this final section, we make a number of remarks about the ideas and results of the earlier sections as well as focus on some open questions and conjectures that we hope to address. We will also mention some results in a forthcoming paper with F. Gesztesy [10] that will study the objects of this paper.

1. Our reconstruction procedure is one-sided, as it must be since $m(z, x)$ is a function of $q$ on $[x, b]$ and totally independent of $q$ on $[0, x]$. The one-sidedness comes from the fact that the differential equation for $A$ begins $\frac{\partial A}{\partial x} = \frac{\partial A}{\partial \alpha}$, not $\frac{\partial A}{\partial x} = -\frac{\partial A}{\partial \alpha}$. If one took an $m_-$ function defined from the left of an interval and normalized so the Riccati equation (1.8) still holds, then $m_-(-\kappa^2)$ has leading asymptotics $+\kappa$ rather than $-\kappa$, and that leads precisely to leading asymptotics $\frac{\partial A}{\partial x} = -\frac{\partial A}{\partial \alpha} + \cdots$ consistent with the one-sidedness in the other direction.

2. We owe to Gel′fand [6] the remark that our basic results extend easily to matrix valued $q$'s (and thus to some higher-order systems). One defines $u$ as a matrix and $m(z) = u'(0, z)u(0, z)^{-1}$, in which case $m$ obeys the matrix equation

$$m' = q - z - m^2.$$

$A$ is matrix-valued. Everything goes through without significant changes.

3. One can ask about the relation of our $A$-function to the kernel $K$ of Gel'fand-Levitan (see 13]). In terms of the Gel'fand-Levitan kernel $K(x, y)$ (defined if $|y| \leq x$), one can define new kernels $K_C, K_S$ defined on $0 \leq y \leq x$ (and built out of $K(x, \pm y)$) so that there are solutions $C, S$ of $-u'' + qu = -\kappa^2 u$ of the form,

$$C(x, \kappa) = \cosh(\kappa, x) + \int_0^x K_C(x, y) \cosh(\kappa y) \, dy$$

$$S(x, \kappa) = \frac{\sinh(\kappa x)}{\kappa} + \int_0^x K_S(x, y) \frac{\sinh(\kappa y)}{\kappa} \, dy.$$

$C, S$ are normalized so that $u_+ = C + m_+ S$, and so defining $u_+$ by the boundary condition at $b$, one gets

(8.1) $$m_+(\kappa) = \frac{hC(b, \kappa) - C'(b, \kappa)}{S'(b, \kappa) - hS(b, \kappa)}.$$

Now,

$$2e^{-\kappa b}(-C' + hC) = -\kappa + h + \kappa \int_0^b B_1(\alpha) e^{-2\alpha\kappa} \, d\alpha$$

$$= -\kappa \left(1 + \int_0^b B(\alpha) e^{-2\alpha\kappa} \, d\alpha \tilde{O}(e^{-2b\kappa})\right)$$

for suitable $B$ defined in terms of $K$ and $h$ and its derivatives. Similarly,

$$2e^{-\kappa b}(S' - hS) = 1 + \int_0^b D(\alpha) e^{-2\alpha\kappa} \, d\alpha + \tilde{O}(e^{-2b\kappa}).$$

By Theorem A.2.3, $(1 + \int_0^b D(\alpha) e^{-2\alpha\kappa} \, d\alpha)^{-1}$ has the form $1 + \int_0^b E(\alpha) e^{-2\alpha\kappa} + \tilde{O}(e^{-2b\kappa})$ and so we can deduce a representation

$$m_+(\kappa) = -\kappa \left(1 + \int_0^b F(\alpha) e^{-2\alpha\kappa} \, d\alpha + \tilde{O}(e^{-2b\kappa})\right).$$

More careful analysis shows that $F(0) = 0$ and $F$ can be differentiated so that $m_+(\kappa) = -\kappa - \int_0^b A(\alpha) e^{-2\alpha\kappa} \, d\alpha + \tilde{O}(\dots)$.

That is, one can discover the existence of our basic representation from the Gel'fand-Levitan representation; indeed, we first found it this way. Because of the need to invert $(1 + \int_0^b D(\alpha) e^{-2\alpha\kappa} \, d\alpha)$, the formula relating $A$ to $K$ is extremely complicated. Subsequent to the preparation of this paper, Gesztesy and I [10] found a simple relation between $A$ and the second Gel'fand-Levitan kernel, $L$, related to $K$ by $1 + L = (1 + K)^{-1}$.

4. The discrete analog of $A$ is just the Taylor coefficients of the discrete $m$-function at infinity. There is, of course, a necessary and sufficient condition for such a Taylor series to come from a discrete Jacobi matrix $m$-function. For these Taylor coefficients are precisely the moments of the spectral measure, and there are a set of positivity conditions such moments have to obey. This suggests that $A$ must obey some kind of positivity conditions. What are they? Is there perhaps a beautiful theorem that the differential equation obeyed by the $A$-function has a solution with a given initial condition if and only if these positivity conditions are obeyed? Subsequent to the preparation of this paper, Gesztesy and I [10] found a simple relation between $A$ and the spectral measure, which is the analog of the Taylor coefficient,

$$A(\alpha) = -2 \int_0^\infty \frac{d\rho(\lambda)}{\lambda^{1/2}} \sin(2\alpha\sqrt{\lambda}),$$

where the divergent integral has to be interpreted as an Abelian limit.

5. The sequence of $\delta$ and $\delta'$ singularities that occur when $b < \infty$ must be intimately related to the distribution of eigenvalues of the associated $H$ via some analog of the Poisson summation formula.

6. There must be an analog of the approach of this paper to inverse scattering theory. Find it!

7. In [10], Gesztesy and I will compute the $A$-function in case $q(x) = -\gamma$ for some $\gamma > 0$. Then

$$A(\alpha) = \frac{\sqrt{\gamma}}{\alpha} I_1(2\alpha\sqrt{\gamma}),$$

where $I_1$ is the standard Bessel function denoted by $I_1(\,\cdot\,)$. Since

$$I_1(z) = \tfrac{1}{2} z \sum_{k=0}^\infty \frac{(\tfrac{1}{4}z^2)^k}{k!(k+1)!},$$

the $\frac{1}{n!}$ bounds in (2.7) are not good as $n \to \infty$ if $q$ is bounded. This is discussed further in [10].

## Appendix 1: Localization of asymptotics

Our goal in this appendix is to prove one direction of Theorem 1.2, viz.:

THEOREM A.1.1. *If $(q_1, b_1, h_1), (q_2, b_2, h_2)$ are two potentials and $a < \min(b_1, b_2)$ and if*

(A.1.1) $\qquad\qquad q_1(x) = q_2(x) \qquad \text{on } (0, a),$

*then as $\kappa \to \infty$,*

(A.1.2) $\qquad\qquad m_1(-\kappa^2) - m_2(-\kappa^2) = \tilde{O}(e^{-2\kappa a}).$

While we know of no explicit reference for this form of the result, the closely related Green's function bounds have long been in the air, going back at least to ideas of Donoghue, Kac, and McKean over thirty years ago. A basic role in our proof will be played by the Neumann analog of the Dirichlet relation (2.2). Explicitly, if $G^D(x,y;z,q)$ and $G^N(x,y;z,q)$ are the integral kernels of $(H-z)^{-1}$ with $H = -\frac{d^2}{dx^2} + q(x)$ on $L^2(0,\infty)$ with $u(0) = 0$ (Dirichlet) and $u'(0) = 0$ (Neumann) boundary conditions, respectively, then

$$m(z) = \lim_{\substack{x<y \\ y\downarrow 0}} \frac{\partial^2 G^D}{\partial x \partial y} \tag{A.1.3}$$

and

$$m(z) = [-G^N(0,0;z,q)]^{-1}. \tag{A.1.4}$$

To see this, let $u$ be the solution $L^2$ at $\infty$ (or which obeys the boundary condition at $b$) and let $\tilde{w}$ obey $-\tilde{w}'' + q\tilde{w} = z\tilde{w}$ with $\tilde{w}(0) = 1$, $\tilde{w}'(0) = 0$ boundary conditions. Then

$$G^N(x,y;z,q) = -\tilde{w}(\min(x,y)) \frac{u(\max(x,y))}{u'(0)}, \tag{A.1.5}$$

from which (A.1.4) is immediate.

We will begin the proof of Theorem A.1.1 by considering the case where $b_1 = b_2 = \infty$.

PROPOSITION A.1.2. *Let $q_1, q_2$ be defined on $(0,\infty)$ and obey (1.4)/(1.5). Then*

$$G^N(0,0;-\kappa^2,q_i) = \kappa^{-1} + o(\kappa^{-1}) \tag{A.1.6}$$

*and if (A.1.1) holds, then*

$$G^N(0,0;-\kappa^2,q_1) - G^N(0,0;-\kappa^2,q_2) = \tilde{O}(e^{-2\kappa a}). \tag{A.1.7}$$

*Remark.* (A.1.4), (A.1.6), and (A.1.7) imply (A.1.2) in this case.

*Proof.* Let $P(x,y;t,q)$ be the integral kernel of $e^{-tH}$ on $L^2(\mathbb{R}, dx)$ where $H = -\frac{d^2}{dx^2} + q(|x|)$. The method of images implies that for $x, y \geq 0$,

$$G^N(x,y;-\kappa^2,q) = \int_0^\infty [P(x,y;t,q) + P(x,-y;t,q)]e^{-\kappa^2 t}\, dt. \tag{A.1.8}$$

Simple path integral estimates (see [16]) imply that

$$P(0,0;t,q) = (4\pi t)^{-1/2}[1 + o(1)] \qquad \text{as } t \downarrow 0 \tag{A.1.9}$$

and if (A.1.1) holds, then for any $\varepsilon > 0$, there exists $C_\varepsilon > 0$ (depending only on the $\beta_2$ for $q_1, q_2$), so that

$$|P(0,0;t,q_1) - P(0,0;t,q_2)| \leq C_\varepsilon \exp(-(1-\varepsilon)a^2/t). \tag{A.1.10}$$

(A.1.9) implies (A.1.6) since $\int_0^\infty 2(4\pi t)^{-1/2} e^{-\kappa^2 t}\, dt = \kappa^{-1} \int_0^\infty (\pi t)^{-1/2} e^\kappa\, dt = \kappa^{-1}$.

To obtain (A.1.7), we use (A.1.8), (A.1.10), and
$$|P(0,0;t,q_j)| \leq C_1 e^{Dt}$$
since
$$\int_0^1 e^{a^2/t} e^{-\kappa^2 t}\, dt = e^{-2\kappa a} \int_0^1 e^{-(a^{-1/2} - \kappa t^{1/2})^2}\, dt$$
$$= O(e^{-2\kappa a}).$$
□

Next, we consider a situation where $b < \infty$, $q$ is given in $L^1(0,b)$, and $h$ is 0 or $\infty$. Define $\tilde{q}$ on $\mathbb{R}$ by requiring that

$$\tilde{q}(x + 2mb) = \tilde{q}(x) \quad m = 0, \pm 1, \pm 2, \ldots, \text{ all } x \in \mathbb{R}$$
$$\tilde{q}(-x) = \tilde{q}(x) \quad \text{all } x \in \mathbb{R}$$
$$\tilde{q}(x) = q(x) \quad x \in [0,b]$$

which uniquely defines $\tilde{q}$ (since each orbit $\{\pm x + 2mb\}$ contains one point in $[0,b]$). Let $G^{(N,N)}$ and $G^{(N,D)}$ be the Green's functions of $-\frac{d^2}{dx^2} + q(x)$ on $L^2(0,b)$ with $u'(0) = 0$ boundary conditions at zero and $u'(b) = 0$ ($(N,N)$ case) or $u(b) = 0$ ($(N,D)$ case) boundary conditions at $b$. Let $\tilde{G}$ be the Green's function for $-\frac{d^2}{dx^2} + \tilde{q}$ on $L^2(\mathbb{R})$. Let $P$ be the corresponding integral kernels for $e^{-tH}$.

By the method of images for $x, y \in [0,b]$:

(A.1.11) $$G^{(N,N)}(x,y;-\kappa^2) = \sum_{m=-\infty}^{\infty} \tilde{G}(x, i_m(y); -\kappa^2)$$

(A.1.12) $$G^{(N,D)}(x,y;-\kappa^2) = \sum_{m=-\infty}^{\infty} \sigma_m \tilde{G}(x, i_m(y); -\kappa^2),$$

where
$$i_m(y) = y + mb \quad m = 0, \pm 2, \pm \cdots$$
$$= -y + mb + b \quad m = \pm 1, \pm 3, \pm \cdots$$
$$\sigma_m = -1 \quad m = 1, 2, 5, 6, 9, 10, \ldots, -2, -3, -6, -7, \ldots$$
$$= 1 \quad \text{otherwise}$$

(i.e., $\sigma_m = -1$, if and only if $m = 1, 2 \mod 4$).

By a simple path integral (or other) estimate on $P$ and Laplace transform, we have

(A.1.13) $$|\tilde{G}(x,y;-\kappa^2)| \leq C_\varepsilon e^{-\kappa|x-y|(1-\varepsilon)}$$

for any $\varepsilon > 0$ and $\kappa$ sufficiently large. Since the images of 0 are $\pm 2b, \pm 4b, \ldots$, (A.1.11) and (A.1.2) imply

PROPOSITION A.1.3.

(A.1.14) $$|G^{(N,N)}(0,0;-\kappa^2) - \tilde{G}(0,0;-\kappa^2)| = \tilde{O}(e^{-2b\kappa})$$

and similarly for $|G^{(N,D)}(0,0;-\kappa^2) - \tilde{G}(0,0;-\kappa^2)|$.

*Remark.* (A.1.14) and (A.1.6) imply (A.1.2) for the pairs $q_1 = \tilde{q}$, $b_1 = \infty$ and $q_2 = q$, $b_2 = b$, and $h_2 = 0$ or $\infty$.

Finally, we compare $b < \infty$ fixed for any two finite values of $h$:

PROPOSITION A.1.4. *Let $q \in L^1(0,b)$. For $h < \infty$, let $G^h$ be the integral kernel for $(-\frac{d^2}{dx^2} + q - z)^{-1}$ with boundary conditions $u(0) = 0$ and $u'(b) + hu(b) = 0$. Then*

(A.1.15) $$|G^h(0,0;-\kappa^2) - G^{h=0}(0,0;-\kappa^2)| = \tilde{O}(e^{-2b\kappa}).$$

*Proof.* Let $H$ be the $h=0$ operator and $H_h$ the operator for $h < \infty$. By the analysis of rank one perturbations (see, e.g., [17]),

$$H_h = H + h(\delta_b, \cdot)\delta_b,$$

where $\delta_b \in \mathcal{H}_{-1}(H)$ is the function $(\delta_b, g) = g(b)$.

Again, by the theory of rank one perturbations [17], let $F(z,h) = G^h(b,b;z)$. Then

$$F(z,h) = \frac{F(z,0)}{1 + hF(z,0)}$$

and

(A.1.16)
$$G^h(0,0;z) - G^{h=0}(0,0;z) = -hG^{h=0}(0,b;z)G^{h=0}(b,0;z)[1 - hF(z,h)]$$
$$= -hG^{h=0}(0,b;z)G^{h=0}(b,0;z)[1 + hF(z,0)]^{-1}.$$

Now $F(-\kappa^2, 0) = \kappa^{-1} + o(\kappa^{-1})$ (this is essentially (A.1.6)) while (A.1.11) and (A.1.13) imply that

(A.1.17) $$G^{h=0}(0,b;z) = \tilde{O}(e^{-\kappa b}).$$

(A.1.16) and (A.1.17) imply (A.1.15). $\square$

Transitivity and Propositions A.1.2–A.2.4 imply Theorem A.1.1.

We close the appendix with two remarks:

1. Do not confuse the Laplace transform in (1.24) (which is in $2\kappa$) with that in (A.1.8) (which is $\kappa^2$).

2. We used path integrals above. As long as $q(x) = O(e^{b|x|})$ for some $b < \infty$, one can instead use more elementary Green's function estimates.

### Appendix 2: Some results on Laplace transforms

In this paper, I need some elementary facts about Laplace transforms. While I am sure that these facts must be in the literature, I was unable to locate them in the precise form needed, so I will give the simple proofs below.

LEMMA A.2.1. *Let $f \in L^1(0, a)$. Suppose that $g(z) \equiv \int_0^a f(y) e^{-zy} \, dy$ obeys*

$$g(x) = \tilde{O}(e^{-ax}) \tag{A.2.1}$$

*as $x \to \infty$. Then $f \equiv 0$.*

*Proof.* Suppose first that $f$ is real-valued. $g(z)$ is an entire function which obeys

$$|g(z)| \leq \|f\|_1 e^{a\operatorname{Re}_-(z)},$$

where $\operatorname{Re}_-(z)$ is the negative part of $\operatorname{Re} z$. Moreover, along the real axis, $g$ obeys (A.2.1). Because of this,

$$h(w) = \int_0^\infty g(x) e^{iwx} \, dx$$

is an analytic function of $w$ in the region $\operatorname{Im} w > -a$. Now for $r > 0$:

$$\begin{aligned}
h(ir) &= \int_0^\infty g(x) e^{-rx} \, dx \\
&= \int_0^\infty \left( \int_0^a f(y) e^{-x(y+r)} \, dy \right) dx \\
&= \int_0^a \frac{f(y)}{y + r} \, dy,
\end{aligned} \tag{A.2.2}$$

where the interchange of integration variables is easy to justify. (A.2.2) implies that

$$h(w) = \int_0^a \frac{f(y)}{y - iw} \, dy \tag{A.2.3}$$

holds for $w$ with $\operatorname{Im} w > 0$ and then allows analytic continuation into the region $\mathbb{C}\setminus\{is \mid s < 0\}$. (A.2.3) and the reality of $f$ implies that for almost every $r \in (0, a)$, $f(r) = \lim_{\varepsilon \downarrow 0} \frac{1}{2\pi i}[h(\varepsilon - ir) - h(-\varepsilon - ir)]$, so the analyticity of $h$ in $\operatorname{Im} w > -a$ implies that $f \equiv 0$. For general complex valued $f$, consider the real and imaginary parts separately. $\square$

An immediate consequence of this is the uniqueness of inverse Laplace transforms.

THEOREM A.2.2. *Suppose that $f, g \in L^1(0, a)$ and for some $b \leq a$, $\int_0^a f(y) e^{-xy} \, dy - \int_0^a g(y) e^{-xy} \, dy = \tilde{O}(e^{-bx})$. Then $f \equiv g$ on $[0, b)$.*

The other fact we need is that the set of Laplace transforms has a number of closure properties. Let $\mathcal{L}_a$ be the set of functions, $f$, analytic in some region $\{z \mid |\operatorname{Arg}(z)| < \varepsilon\} \equiv R_\varepsilon$ obeying

$$f(z) = 1 + \int_0^a g(\alpha) e^{-\alpha z} \, d\alpha + \tilde{O}(e^{-a \operatorname{Re} z})$$

in that region for some $g \in L^1(0, a)$. Denote $g$ by $\mathcal{I}(f)$.

THEOREM A.2.3. *If $f, h \in \mathcal{L}_a$ so are $fh$, $f + h - 1$, and $f^{-1}$.*

*Proof.* $f + h - 1$ is trivial. $fh$ is elementary; indeed,

$$\mathcal{I}(fh)(\alpha) = \mathcal{I}(f)(\alpha) + \mathcal{I}(h)(\alpha) + \int_0^\alpha \mathcal{I}(f)(\beta) \mathcal{I}(h)(\alpha - \beta) \, d\beta.$$

For the inverse, we start by seeking $k$ obeying (where $g = \mathcal{I}(f)$)

$$g(\alpha) + k(\alpha) + \int_0^\alpha d\beta \, k(\beta) g(\alpha - \beta) = 0.$$

This Volterra equation always has a solution (by iteration). Let $h(z) = 1 + \int_0^a k(\alpha) e^{-\alpha z} \, d\alpha$. Then

$$fh = 1 + \tilde{O}(e^{-a \operatorname{Re}(z)})$$

and so

$$\begin{aligned} f^{-1} &= h(1 + \tilde{O}(e^{-a \operatorname{Re}(z)}))^{-1} \\ &= h + \tilde{O}(e^{-a \operatorname{Re}(z)}) \end{aligned}$$

as required. $\square$

*Notes added in proof.*

1. For the case of short-range potentials, a representation of the form (2.3) was obtained by A. Ramm in the paper, "Recovery of the potential from the $I$-function," C. R. Math. Rep. Acad. Sci. Canada **IX** (1987), 177–182.

2. Recently, F. Gesztesy and the author obtained an alternate and simpler proof of Theorem 1.2 in the paper, "On local Borg-Marchenko uniqueness results," which will appear in Commun. Math. Physics.

California Institute of Technology, Pasadena, CA
*E-mail address*: bsimon@caltech.edu